\documentclass[a4paper,
10pt
]{article}

\usepackage[a4paper,left=3cm,right=3cm,top=2.5cm,bottom=2.5cm]{geometry}
\usepackage{hyperref}
\usepackage{subfig}
\usepackage{pgfplots}
\usepackage{pgfplotstable}
\usepackage{mathtools}
\usepackage{multicol}
\usepackage{comment}
\usepackage{booktabs}
\pgfplotsset{compat=1.5}
\usepackage{amssymb}
\usepackage{url}
\usepackage{bm}

\usepackage{pdfsync}
\usepackage{float}
\usepackage{tabularx}
\usepackage{enumerate}
\usepackage{array}
\usepackage{xspace}
\usepackage{tikz}
\usepackage{tikz-cd}
\usepackage{tikzsymbols}
\usetikzlibrary{calc,trees,positioning,arrows,chains,shapes.geometric,%
    decorations.pathreplacing,decorations.pathmorphing,shapes,%
    matrix,shapes.symbols, decorations.markings, patterns,fit}

\usepackage{siunitx}

\usepackage{authblk}

\usepackage{mathrsfs}
\usepackage{color, colortbl}
\usepackage{multirow}

\usepackage{amsthm}
\usepackage{amsmath}
\usepackage{stmaryrd}
\usepackage[ruled,vlined,linesnumbered]{algorithm2e}
\usepackage{graphicx}
\usepackage{booktabs}
\usepackage[noabbrev]{cleveref}

\newtheorem{theorem}{Theorem}

\theoremstyle{remark}

\newcommand{\review}[1]{{\color{black}#1}}

\newcommand{\R}{\ensuremath{\mathbb{R}}}

\newcolumntype{C}[1]{>{\centering\arraybackslash}m{#1}}

\definecolor{Gray}{gray}{0.9}

\newcommand{\maria}[1]{\textbf{\color{teal}{#1}}\color{black}}
\newcommand{\dual}{^\ast}
\newcommand{\Cal}{\mathcal}
\newcommand{\norm}[1]{\lVert #1 \rVert}
\newcommand{\half}{\frac{1}{2}}
\newcommand{\alf}{\frac{\alpha}{2}}

\newcommand{\control}{\mathcal U}
\newcommand{\la}{\langle}
\newcommand{\ra}{\rangle}
\newcommand{\ocp}{OCP}

\newcommand{\intTime}[1]{\int_0^T  #1 \; dt}

\newcommand{\dt}[1]{{ #1}{_t}}
\newcommand{\dn}[1]{\frac{\partial #1}{\partial n}}

\newcommand{\goesto}[3]{: #1 \times #2 \rightarrow #3}
\newcommand{\Lg}{\mathscr{L}}

\begin{document}

\title{A data-driven partitioned approach for the resolution of
  time-dependent optimal control problems with dynamic mode
  decomposition} 

\author[]{Eleonora~Donadini\footnote{edonadini@sissa.it}}
\author[]{Maria~Strazzullo\footnote{maria.strazzullo@sissa.it}}
\author[]{Marco~Tezzele\footnote{marco.tezzele@sissa.it}}
\author[]{Gianluigi~Rozza\footnote{gianluigi.rozza@sissa.it}}

\affil{Mathematics Area, mathLab, SISSA, via Bonomea 265, I-34136
  Trieste, Italy}

\maketitle

\begin{abstract}
This work recasts time-dependent
optimal control problems governed by partial differential equations in
a Dynamic Mode Decomposition with control framework. Indeed, since the
numerical solution of such problems requires a lot of computational
effort, we rely on this specific data-driven
technique, using both solution and desired state measurements to
extract the underlying system dynamics. Thus, after the Dynamic Mode
Decomposition operators construction, we reconstruct and perform
future predictions for all the variables of interest at a lower
computational cost with respect to the standard space-time discretized
models. We test the methodology in terms of relative reconstruction
and prediction errors on a boundary control for a Graetz flow and on a
distributed control with Stokes constraints.  
\end{abstract}

\tableofcontents

\section{Introduction}
\label{sec:intro}
Scientific and industrial contexts need to
represent natural phenomena through mathematical
models. Historically, this role has been played by partial
differential equations (PDEs) and by their numerical simulations. Yet,
the model equations, in some frameworks, are not enough to well
describe the complexity of the natural phenomenon one is dealing
with. To increase the reliability of a proposed PDEs-based model, a
very classical and elegant mathematical tool has been exploited:
optimal control. This technique responds to the need for reducing the
gap between PDEs and collected data. Indeed, the data information,
given by some previous knowledge on the system, is exploited in order
to achieve a \emph{desired configuration}: a convenient profile
similar to the behaviour observed or expected in nature. Namely,
through a PDE-constrained minimization strategy, optimal control
problems (OCPs) are able to steer the solution of the problem at hand
towards a target profile.
\ocp s governed by PDEs are widespread in many scientific
applications: some examples in shape optimization may be found in
\cite{delfour2011shapes,makinen,mohammadi2010applied} or, in fluid
dynamics, in
\cite{dede2007optimal,dede2010reduced,negri2015reduced}
where a parametrized setting is discussed. In the same framework, we can cite biomedical applications
{\cite{Ballarin2017,Carere2021261,Fevola2021,LassilaManzoniQuarteroniRozza2013a,ZakiaMaria,Zakia}}, or environmental ones
\cite{BALLARIN2022307,quarteroni2005numerical,quarteroni2007reduced,Strazzullo1,StrazzulloSWE}. 
From this references, it is clear that OCPs are of indisputable usefulness in many applications. Yet,
they still are very complex to analyse and to simulate, most of all if
time-dependency is taken into consideration, since it requires larger
computational costs.
Despite these difficulties, time-dependent \ocp s have been treated in
many works: here a far from exhaustive list
\cite{HinzeStokes,Iapichino2,leugering2014trends,seymen2014distributed,Stoll,Strazzullo2,StrazzulloRB,StrazzulloSWE}.
In contrast with the aforementioned works, we tackle the study
of the evolution of \ocp s trough the employment of dynamic mode
decomposition with control
(DMDc)~\cite{proctor2016dynamic}, which is a data-driven technique for
system identification in the context of feedback and control.
The idea is to use the desired state from the \ocp s
framework as forcing term within the DMDc. Then, with a partitioned
approach, we use the solution snapshots to construct different DMD
operators. With such operators, given new actuation snapshots, we
perform future predictions for all the variables of 
interest, in a data-driven fashion. This results in great 
advantages in terms of computational time, keeping the accuracy within
a certain threshold, as usual in the reduced order modeling framework~\cite{morhandbook2020}.
The main novelty of this work, thus, is the recasting
of an OCP in a DMDc framework: to the best of our knowledge, indeed,
this is the first time that it is highlighted and investigated. We
propose two test cases: a boundary control governed by a Graetz flow and a
distributed control with Stokes constraints. 

\section{Time-dependent optimal control problems}
\label{sec:ocp}
We provide the continuous formulation for general time dependent OCP. Let us consider the spatial domain $\Omega \subset
\R^d$, with $d=2,3$: here, the analysed physical phenomenon,
described by a linear time-dependent PDE, is taking place in the time
interval $[0,T]$. Furthermore, we denote with $\Gamma_D$ and $\Gamma_N$ two non-overlapping boundary portions where homogeneous Dirichlet and Neumann boundary conditions apply. In the context of constrained optimization, we consider an Hilbert space $Y$ such that $Y \hookrightarrow H \hookrightarrow Y\dual$ for some suitable $H$. Moreover, we define 
$$
\mathcal{Y}_t = \{y \in L^2(0,T; Y) \text{ such that } y_t \in  L^2(0,T; Y\dual)  \text{ with } y(t) = 0  \} \subset \mathcal{Q},
$$
with $\mathcal Q = L^2(0,T; Y)$. The problem \emph{state variable} $y$ is sought in $\mathcal Y_0$. In addition, we need to define
$\mathcal U = L^2(0,T; U)$ as the space for the \emph{control
  variable} $u$, with $U$ another suitable Hilbert space. \\The control
acts on $\Omega_u \subseteq \Omega$ of a portion of its boundary, say
$\Gamma_C \subset \partial \Omega$ with $\Gamma_C \cap \Gamma_N \cap \Gamma_D = \emptyset$ while $\Gamma_D \cup \Gamma_N \cup \Gamma_C = \partial \Omega$. In the first case, we say that the
control problem is \emph{distributed}, while in the second case, we
say that the problem is a \emph{boundary control}. With no distinctions, from now on, we will call $\Omega_u$ or $\Gamma_C$ the \emph{control domain}.
\\  In order to change the classical solution behaviour, we need to define a \emph{controlled system} of the following form, for all $q$ in $\mathcal Q$, considering once again $y_0$ as initial condition
\begin{align}
\label{eq_time_weak}
& \displaystyle \intTime{ \left \la \dt y, q \right \ra_{Y\dual, Y}} \displaystyle + \intTime {a (y, q)} = 
 \intTime {c(u, q)} +
\intTime{ \la G, q \ra_{Y\dual, Y}}, 
\end{align}
with initial condition $y_0$ where, $a \goesto{Y}{Y}{\mathbb R}$ is a coercive and continuous bilinear form and 
$ G \in Y\dual $ represents the forcing and the boundary terms
of the problem at hand, and $c \goesto{U}{Y}{\mathbb R}$ is the $L^2$ product on the {control domain}.
For the sake of notation, we define the controlled equation as 
$\mathcal E: (\mathcal Y \times \mathcal U) \times \mathcal Q  \rightarrow \mathbb R$. Namely, \eqref{eq_time_weak} is verified when $
    \mathcal E((y,u), q) = 0$ for all $ q$ in $\mathcal Q.$ \\
The goal of an OCP is to steer the solution towards a desired profile $y_\text{d} \in {\mathcal Z} \supseteq \mathcal Y$ by minimizing a convex \emph{cost functional} $J: \mathcal Y_0 \times \mathcal U \rightarrow \mathbb R$ over $\mathcal Y_0 \times \mathcal U$ such that \eqref{eq_time_weak} is verified. The functional $J(y,u)$ must be defined case by case. For the well-posedness of the problem, the interested reader may refer to \cite{hinze2008optimization}. The minimization problem can be recast in an unconstrained fashion exploiting a Lagrangian approach \cite{hinze2008optimization,troltzsch2010optimal}. Thus, let us consider an arbitrary \emph{adjoint variable} $z \in {\mathcal Y_T} \subset \mathcal Q$
and define the following Lagrangian functional
$$\Lg (y,u,z) = J(y,u) + \Cal E((y,u), z)$$. It is well known in literature~\cite{hinze2008optimization}, that the
aforementioned minimization of problem is equivalent to the
following system: find $(y,u, z) \in \Cal Y_0 \times \Cal U \times \Cal Y_T$ such that 
\begin{equation}
\label{eq:optimality_system}
\begin{cases}
D_y\Lg(y, u, z)[\omega] = 0 & \forall \omega \in \Cal Q \quad \text{(adjoint equation)},\\
D_u\Lg(y, u, z)[\kappa] = 0 & \forall \kappa \in \Cal U \quad \text{(optimality equation)},\\
D_z\Lg(y, u, z) [\zeta]= 0 & \forall \zeta \in \Cal Q \quad \text{(state equation)}.\\
\end{cases}
\end{equation}
In this context, the terms $D_y, D_u$ and $D_z$ will denote the
differentiation with respect to the state, the control and the adjoint
variables.
We underline that, over the control domain, the optimality equation in strong form reads:
\begin{equation}
\label{eq:alpha_and_p}
\alpha u - z = 0,
\end{equation}
where the role of $\alpha$ will be clarified in the next sections.
We decided to solve~\eqref{eq:optimality_system} using a \emph{space-time discretization}. The space-time techniques have been successfully applied in many contexts, from parabolic equations \cite{Glas2017,urban2012new,yano2014space,yano2014space1} to time-dependent PDE-constrained optimization~\cite{HinzeStokes,HinzeNS,Strazzullo2,StrazzulloSWE,StrazzulloRB}. For the sake of brevity, we are not presenting the details of such an approach. The main idea is to treat the three-equations system \eqref{eq:optimality_system} as a steady one, where all the time instances are sought all-at-once by means of a direct solver. The reader may refer to the previous cited bibliography for an insight to the method. \\For our purposes, it is enough to define a time discretization over  $[0,T]$ divided in $N_t$ sub-intervals. Namely, we consider $\Delta t > 0$ and the time instance $t_k = k\Delta t$ for $k = 0, \dots, N_t$. 
Let us focus on the state, adjoint and control variables evaluated at $t_k$. At each $t_k$, the variables can be expressed through the corresponding spatial bases, in our case Finite Element (FE) bases.
From now on, for the sake of notation, $y_k$ will denote the column vector of FE coefficients of the FE expansion.  The same argument applies to control and adjoint variables at $t_k$, denoted by $u_k$ and $p_k$. We remark that state and adjoint have been discretized with the same discretized function space \cite{Strazzullo2} to guarantee the well-posedness of the problem \eqref{eq:optimality_system}. To simulate the solution of the OCP we solve a system of dimension $N_t (2\mathcal N_y + \mathcal N_u) \times N_t (2\mathcal N_y + \mathcal N_u)$ through a direct solver. This may lead to time consuming simulations to understand the field dynamic. To lighten these computational costs, one may use multigrid approaches combined with proper preconditioners, see e.g.\ \cite{Benzi,schoberl2007symmetric,Stoll} and the references therein. However, these iterative methods might only partially solve this issue. DMD-based techniques respond to the need of a fast prediction tool to avoid a complete simulation over the time interval $[0,T]$. 

\section{Dynamic mode decomposition with control}
\label{sec:dmd}

Dynamic mode decomposition (DMD) is a powerful method to identify and
approximate dynamical systems using only few spatiotemporal coherent
structures~\cite{brunton2019data,kutz2016dynamic}. It is ideally
suited for time-dependent problems, and its data-driven nature makes
it very versatile, also in industrial contexts~\cite{demo2018shape,tezzele2018ecmi,tezzele2020enhancing}.

When dealing with actuated systems, unfortunately, DMD fails to properly reconstruct
the underlying dynamics since it is incapable of incorporating the
contribution of the forcing term. To this end, dynamic mode decomposition with
control (DMDc)~\cite{proctor2016dynamic} has been developed to
overcome this issue. By including the actuation snapshots in the analysis it is able
to provide reduced order representations for input-output systems. The
DMDc method quantifies the effect of control inputs on the state of
the system and computes the underlying dynamics without being
confounded by the effect of external control. 
Recently, an extension for quantum control problems, called bilinear
DMD, has been proposed in~\cite{goldschmidt2021bilinear}. A local
version of DMDc for predictive control of hydraulic fracturing can be
found in~\cite{narasingam2017development}, while for compressive
system identification see~\cite{bai2020dynamic}.

Let us denote the snapshot representing the state of a
system at the $i$-th time instant with $x_i \in \mathbb{R}^{\mathcal{N}}$, where
$\mathcal{N}$ represents the number of degrees of freedom of our system. We
collect a set of $N_t$ equispaced snapshots $\{x_i\}_{i=1}^{N_t}$ and we
arrange them by column in two matrices, $X$ and $X^\prime$, where
$X^\prime$ is the time-shifted version of $X$. We also
collect the input control snapshots $\{\eta_i\}_{i=1}^{N_t-1}$, with $\eta_i
\in \mathbb{R}^{\mathcal{L}}$, in the matrix $\Upsilon$, obtaining the following matrices:
\begin{equation*}
  X=\left[
    \begin{array}{cccc}
      | & | &  & |\\
      x_1 & x_2 & \dots & x_{N_t-1}\\
      | & | &  & |
    \end{array}
  \right],\quad
  X^\prime=\left[
    \begin{array}{cccc}
      | & | &  & |\\
      x_2 & x_3 & \dots & x_{N_t}\\
      | & | &  & |
    \end{array}
  \right],\quad
  \Upsilon=\left[
    \begin{array}{cccc}
      | & | &  &| \\
      \eta_1 & \eta_2 & \dots & \eta_{N_t-1} \\
      | & | &  &|
    \end{array}
  \right].
\end{equation*}
We remark that the snapshots can represent data coming from experiments, from
simulations, or even sensors and acquired in real-time. 

DMDc is a regression-based approach to system identification that is
able to disambiguate the intrinsic dynamics, described by the matrix
$A$, and the effects of control, described by the matrix $B$, as in
the following
\begin{equation}
\label{eq:augmented_matrix_dmd}
X^\prime \approx A X + B \Upsilon =
\left[ \begin{array}{cc} A & B \end{array} \right]
\left[ \begin{array}{c} X \\ \Upsilon \end{array}\right]
:= G \Xi .
\end{equation}
We seek an approximation of the linear mappings $A$ and $B$ using only
the three data matrices. We start by computing the SVD of the matrix
$\Xi$ in~\ref{eq:augmented_matrix_dmd}, which contains both the state and
control snapshot information, as $\Xi \approx U \Sigma V^*$. With the
symbol $^*$ we denote the conjugate transpose, and we truncate the SVD
keeping only the first $r_{\Xi}$ modes. We can thus express the
best-fit matrix $G \in \mathbb{R}^{\mathcal{N} \times (\mathcal{N} + \mathcal{L})}$ as
\begin{equation}
G \approx X^\prime \Xi^{-1} =  X^\prime V \Sigma^{-1} U^*.
\end{equation}
By rewriting $U^* = \left[ \begin{array}{cc} U^*_1 & U^*_2 \end{array}
\right]$, and by computing the SVD of the output matrix $X^\prime
\approx U_{X^\prime} \Sigma_{X^\prime} V_{X^\prime}^*$, we can
compute the reduced order approximation of $A$ and $B$ as
\begin{align}
\tilde{A} &= U_{X^\prime}^* X^\prime V \Sigma^{-1} U^*_1 U_{X^\prime}, \\ 
\tilde{B}  &= U_{X^\prime}^* X^\prime V \Sigma^{-1} U^*_2.
\end{align}
We remark that the truncation rank $r_{X^\prime}$ of the SVD of the
output matrix has to be less or equal to $r_{\Xi}$. With this reduced
order operators we can write
$\tilde{x}_{k+1} = \tilde{A} \tilde{x}_k + \tilde{B} \eta_k$,
where $\tilde{x}_k = U_{X^\prime} x_k$. The dynamic modes of $A$ can
be computed from the eigenvectors $w$ of $\tilde{A}$ as
\begin{equation}
\phi = X^\prime V \Sigma^{-1} U^*_1 U_{X^\prime} w.
\end{equation}

\section{A partitioned approach for time-dependent OCPs}
\label{sec:partitioned}

In this section we are going to explain how to handle \ocp s with DMDc. Our goal is to
characterize the relationship between the current measurements of the
time-dependent OCP $x_k$, the future one $x_{k+1}$, and the current
input $\eta_k$, given by the desired state ${y}_d$ at
each time instant $t_k$ for a finite number of steps $k=1, \dots, N_t$: 
$$x_{k+1}=A x_k + B \eta_k, \quad \forall k = 1, \dots, N_t-1$$.
We remark that the current input for our OCP system is
the desired state ${y}_d$. In fact, it represents the
target state that we want to achieve and hence all the OCP variables
(state, control, and adjoint) depend on it. Due to such a dependency, a
modification in the desired state will cause a change in the other
variables. 
We present a \textit{partitioned approach} in order to identify two
separate dynamical systems: one for the state, and one for the adjoint
variables. The control can then be reconstructed through $\alpha$ and
the adjoint variable by using the linear relation in~\eqref{eq:alpha_and_p}.
For the state variable, we use as measurement matrices $X_{y}\in
\mathbb{R}^{\mathcal{N}_y\times (N_t-1)} $ and $X^\prime_y\in
\mathbb{R}^{\mathcal{N}_y\times (N_t-1)}$ given by: 
\begin{equation}
	X_{y}=\left[\begin{array}{cccc}
		\mid & \mid & & \mid \\
		y_1 & y_2 & \ldots & y_{N_t-1} \\
		\mid & \mid & & \mid
	\end{array}\right], \qquad
	X^\prime_y =\left[\begin{array}{cccc}
		\mid & \mid & & \mid \\
		y_2 & y_3 & \ldots & y_{N_t} \\
		\mid & \mid & & \mid
	\end{array}\right],
\end{equation}
where each row represents the time series of the state measurement for
a particular spatial point given by the space-time discretization. On
the other hand, each column contains the 
coefficients of the FE expansion of the state space-time variable. 
The final model we obtain is:
\begin{equation}
	X^\prime_y = A_y X_y +B_y \Upsilon_d.
\end{equation}

For the adjoint variable we have that $X_{z}\in \mathbb{R}^{\mathcal{N}_y\times (N_t-1)} $, $X^\prime_z\in
\mathbb{R}^{\mathcal{N}_y\times (N_t-1)}$
where each row represents the
time series of the adjoint measurement and each column collects the
coefficients of the FE expansion of the adjoint space-time variable.
We obtain the DMDc model as:
\begin{equation}
	X^\prime_z = A_z X_{z} +B_z\Upsilon_d.
\end{equation}
We emphasize that the input matrix $\Upsilon_d$ is the same for all the models.

\section{Numerical results}
\label{sec:results}
In this section, we present the numerical results to validate the DMDc
approach applied to OCPs. The first example takes into consideration
an unsteady Graetz Flow boundary OCP. A second test case deals with a 
time-dependent OCP governed by Stokes equations. The experiments
validate the performances of DMDc strategy in term of
reconstruction and prediction. Indeed, we define the time-pointwise
relative error $E_k$ as:
\begin{equation}
  \label{eq:relative_error}
E_k =  \frac{\| x_k - \tilde{x}_k \|_2 }{\| x_k \|_2},
\end{equation}
where $x_k$ represents the generic true snapshot at time
$t_k$, and $\tilde{x}_k$ the approximated variable. This error is
shown for the reconstruction analysis. For the prediction, we
average~\eqref{eq:relative_error} over all the time instances in the
test  data set. We vary the size of the training data set and we keep fixed
$20$ as dimension for the test data set, to understand sensitivity of
the method with respect to the data needed for an accurate
prediction. For the computations we used the following libraries:
multiphenics~\cite{multiphenics}, which is an implementation in
FEniCS~\cite{fenics} for block-based systems, and
PyDMD~\cite{demo2018pydmd}. 

\subsection{Boundary OCP governed by a Graetz Flow}
\label{sec:graezt_results}
\begin{figure}[h]
	\centering
	\includegraphics[width=0.99\textwidth]{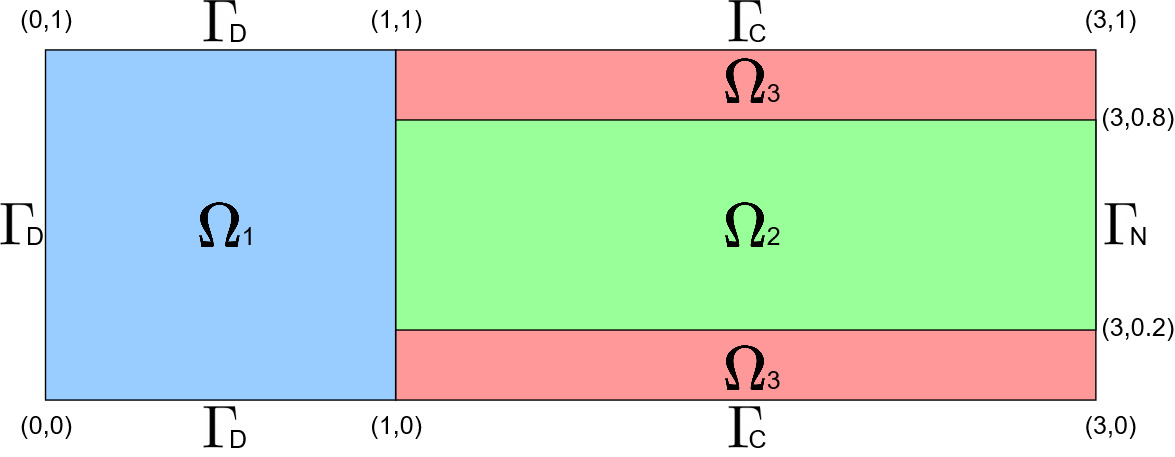}
	\caption{\emph{Graetz Flow}: space domain.}
	\label{fig:rectangularDomain}
\end{figure}

\begin{figure}[htb]
  \centering
  \includegraphics[width=0.49\textwidth]{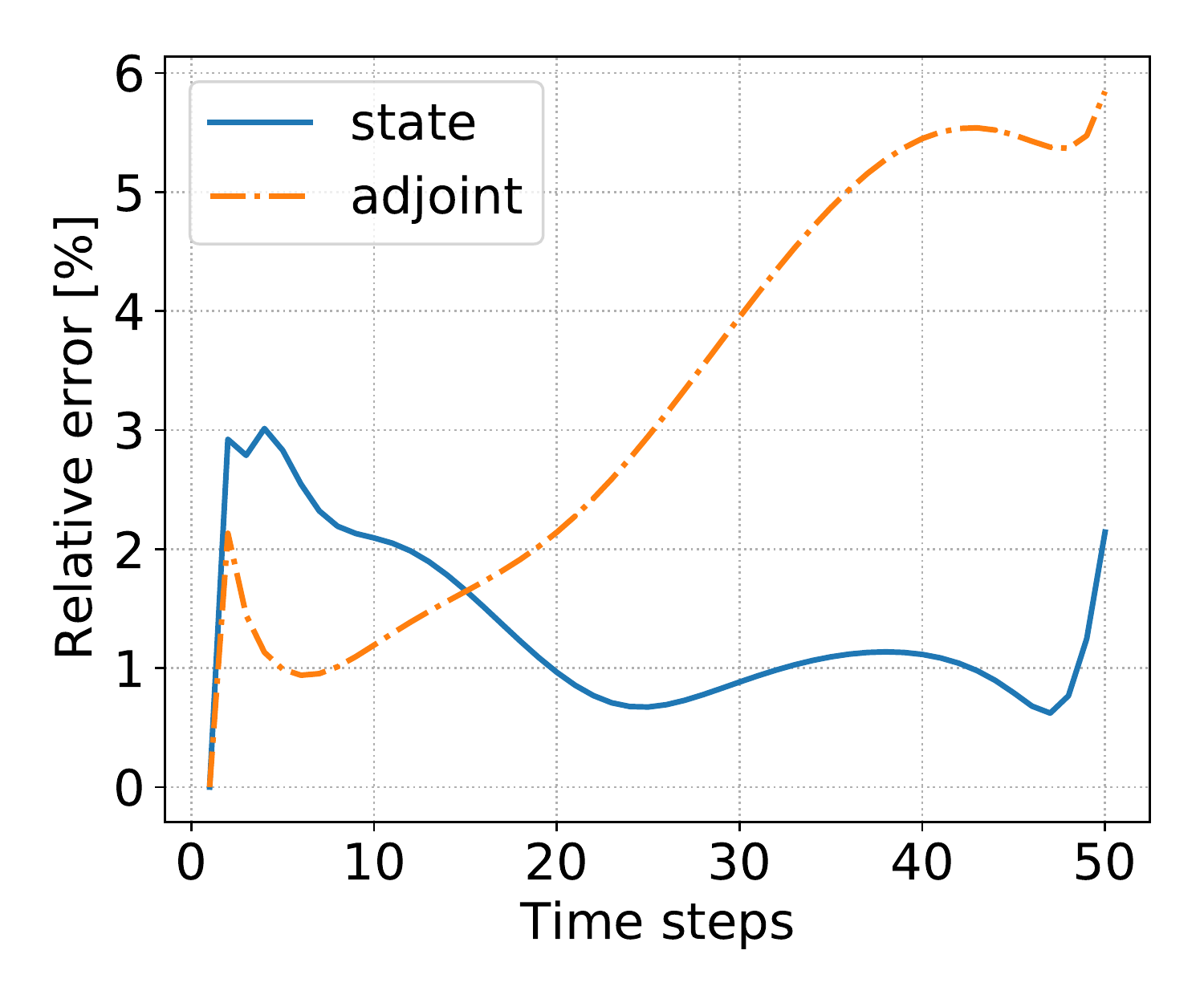} \hfill
  \includegraphics[width=0.49\textwidth]{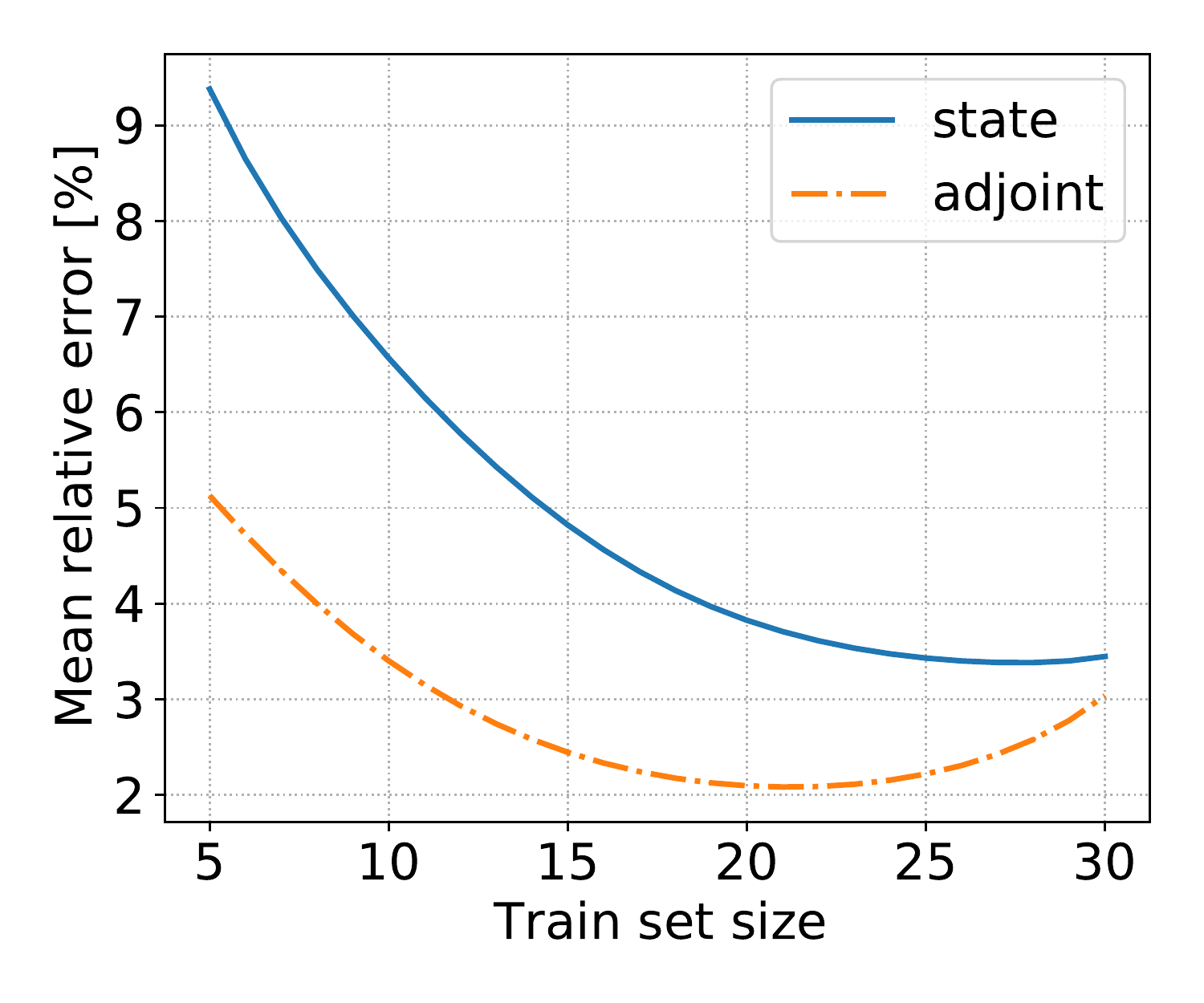}
  \caption{On the left the $L^2$ relative reconstruction error for the Graetz
    case for state and adjoint variables. On the right the $L^2$ mean
    relative prediction error. The average is done over $20$ snapshots in the future.}
  \label{fig:recon_pred_error_boundary}
\end{figure}

Let us consider a boundary control governed by a Graetz Flow in the time interval $[0, T]= [0,1]$ and in the rectangular space domain $\Omega =
[0,3]\times[0,1]\in\mathbb{R}^2$ as depicted in Figure~\ref{fig:rectangularDomain}. 
The control domain is $$\Gamma_C = ([1, 3] \times \{0\}) \cup ([1, 3]
\times \{1\}),$$  while the observation domain is $$\Omega_3 = ([1, 3]\times [0,
0.2]) \cup ([1, 3]\times [0.8,1]).$$ The Neumann and the Dirichlet conditions are applied, respectively, to
$$\Gamma_N = \{3\} \times [0,1] \text{ and }\Gamma_D = \partial \Omega \setminus (\Gamma_C \cup \Gamma_N).$$ We recall that 
$Y=H^1_{\Gamma_D^0}(\Omega)$, i.e.\ the set of functions that are
$H^1$ in the domain but null on $\Gamma_D$, and the control space
$U=L^2(\Gamma_C)$. We recall that we can recover an homogeneous
problem of the form presented in section~\ref{sec:ocp} after a lifting
procedure~\cite{quarteroni2005numerical}. Thus, the Lagrangian
formulation of this OCP reads: find $(y, u) \in \mathcal{Y}_0 \times \mathcal{U}$
which solves: 
\begin{equation}
	\label{eq:adr1}
	\min _{(y, u) \in \mathcal{Y}_0 \times \mathcal{U}} J(y, u) = \frac{1}{2} \int_{0}^T\int_{\Omega_3} (y - y_d)^2 dxdt + \frac{\alpha}{2} \int_{0}^T\int_{\Gamma_C} u^2 dsdt
\end{equation}
constrained to
\begin{equation}
	\label{eq:adr2}
	\begin{cases}
		\dt{y} - \epsilon \Delta y + \beta \cdot \nabla y  = 0 & \text { in }  \Omega \times (0, T),\\
		y =1 & \text{ on } \Gamma_D \times (0, T), \\
		\displaystyle \epsilon \frac{\partial y}{\partial n} =  u & \text{ on }  \Gamma_C \times (0, T),\\
		\displaystyle \epsilon \frac{\partial y}{\partial n} =  0 &  \text{ on }  \Gamma_N \times (0, T),\\
		y = y_0 & \text { in } \Omega\times \{0\},
	\end{cases}
\end{equation}
where $\epsilon=\frac{1}{12}$, the vector filed $\beta$ is defined as
$\left[x_2(1 - x_2), 0\right]$, and $x_2$ as vertical spatial
coordinate. The desired state profile is $y_d = 1 + t$,
$\alpha = 10^{-2}$ is a penalization parameter on the control, and
$y_0$ is a null function verifying the boundary
conditions. For the discretization, we employed
$\mathbb P^1$ elements for all the variables with $\mathcal N_y
=\mathcal N_u = 2304$. Moreover, in terms of time discretization, we
consider $\Delta t = 0.02$, leading to $N_t=50$ time instances for the
time interval $[0,1]$. Thus, the global space-time dimension is
$345600$.  
In the left panel of Figure~\ref{fig:recon_pred_error_boundary} we plotted the relative
reconstruction error $E_k$ for all the variables of interest and for
every $k=1, \dots, 50$. We used $4$ DMD modes for the state and $3$
for the adjoint. The mean relative error for the state is
$1.3\%$, while for the adjoint variable is $3.2\%$.

In Figure~\ref{fig:recon_pred_error_boundary} (right panel) we plotted the $L^2$ mean
relative prediction error over $20$ future states for the state, and
adjoint variable, varying the dimension of the train
dataset. We see that for the state variable we need at least $20$
snapshots in order to have a satisfactory prediction accuracy below
$4\%$, while for the adjoint variable $10$ snapshots are sufficient.

\subsection{OCP governed by time-dependent Stokes Equations}
\begin{figure}[htb]
  \centering
  \includegraphics[width=0.49\textwidth]{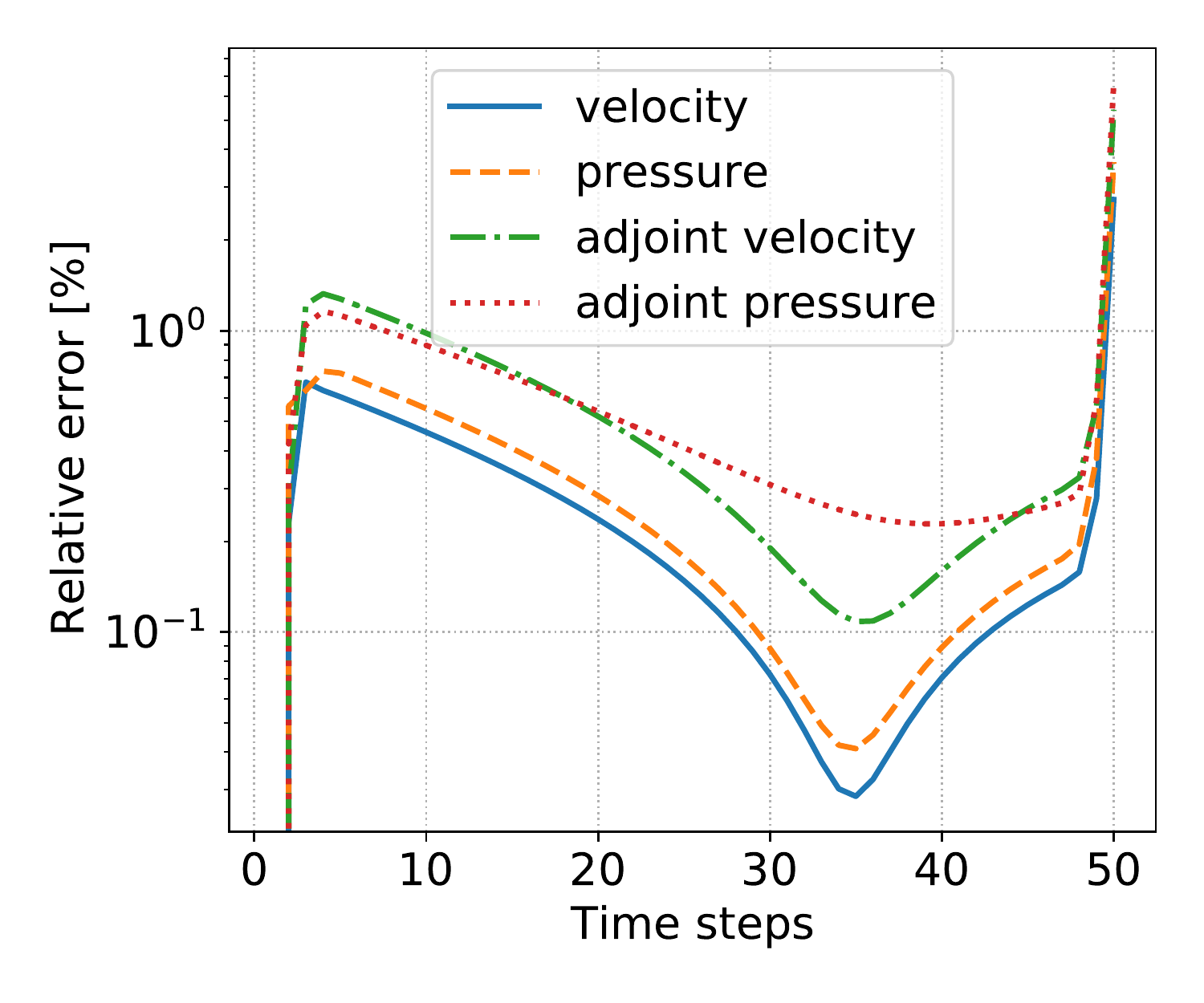}\hfill
  \includegraphics[width=0.49\textwidth]{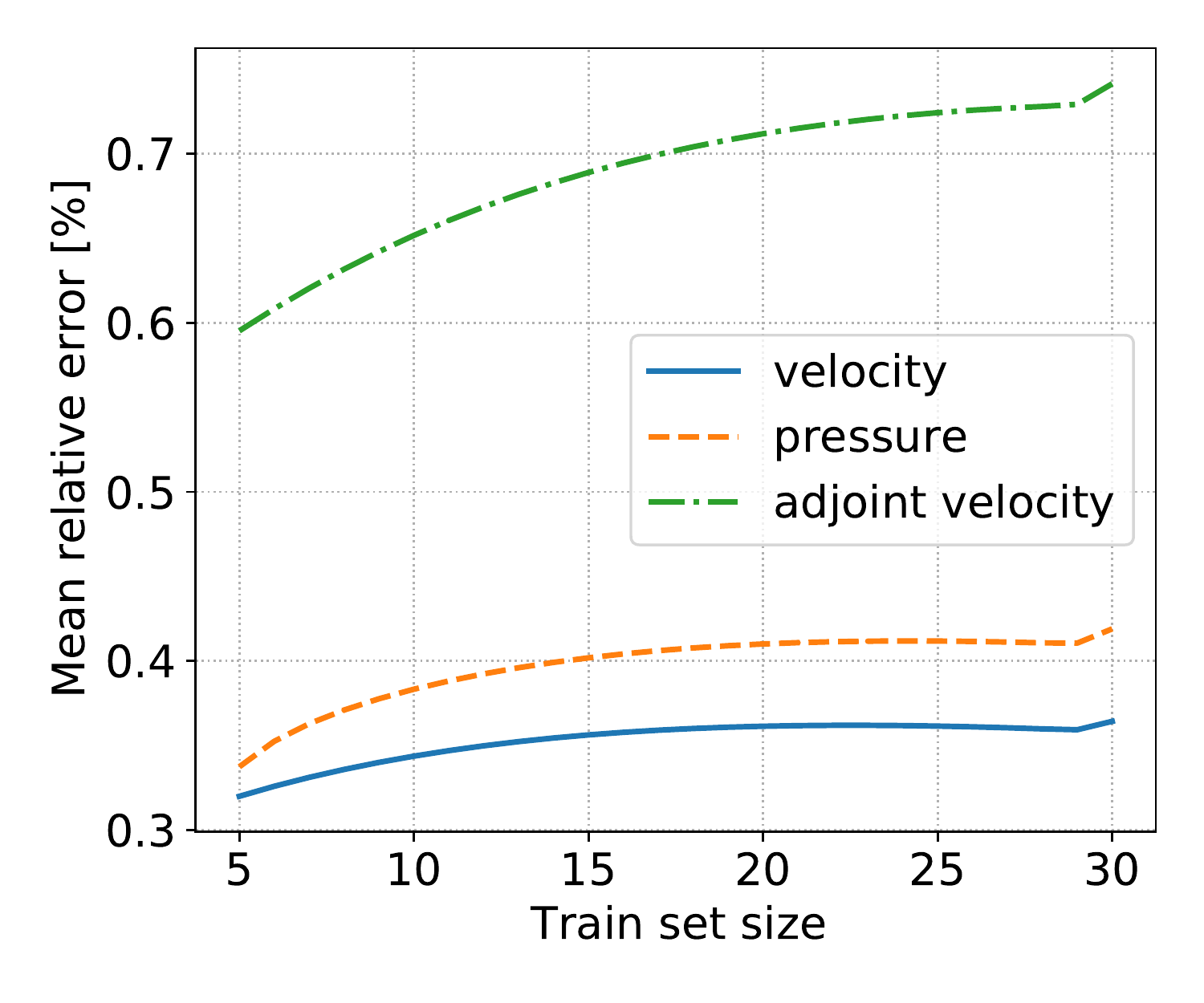}
  \caption{On the left the $L^2$ relative reconstruction error for the Stokes
    case for velocity, pressure, and the adjoint variables. On the
    right the $L^2$ mean relative prediction error. The average is
    done over $20$ snapshots in the future.}
  \label{fig:recon_pred_error_stokes}
\end{figure}

We now deal with a distributed control problem governed by
time-dependent Stokes equation. Let us consider the spatial domain
$\Omega$ as the unit square in $\mathbb{R}^2$ an the time interval
$[0,1]$. For this specific case $\Gamma_D = \partial \Omega$, and we
consider
$$\Cal V_0 \doteq \left \{ v \in L^2(0,T; V) \; \text{such that} \; \displaystyle \dt{v} \in L^2(0,T; V\dual)  \text{ such that } v(0) = 0 \right \} $$
with 
$$ 
V \doteq [H^1_{\Gamma_D}(\Omega)]^2 \text{ and } \Cal P = L^2(0,T;
L^2_0(\Omega)),
$$ where $L^2_0(\Omega)$ is the space of functions that
have null mean over $\Omega$ with $L^2$
regularity\footnote{Numerically, to enforce the condition of null mean
  for the state pressure variable, we employed Lagrange
  multipliers. Namely, the condition is weakly imposed in integral
  form 
\begin{equation}
\label{eq:multiplier1}
\int_{\Omega} p \lambda \; dx = 0, \quad \forall \lambda \in \mathbb R.
\end{equation}
This constraint reflects in the derivation of the adjoint equation and
an extra term appears in the divergence free constraint in weak
form.}. We seek the state-control variable in $\mathcal Y_0 \times
\mathcal U$ where $\mathcal Y_0 = \mathcal V_0 \times \mathcal P$ and
$\mathcal U = [L^2(\Omega)]^2$. 
The specific problem we are dealing with is to find the minimum of 
\begin{equation}
	\min_{((v,p),u) \in \mathcal Y_0 \times \mathcal U}
        J(((v,p), u)):=\frac{1}{2} \int_{0}^{T}
        \int_{\Omega}(v-{v}_d)^{2} d x d t+\frac{\alpha}{2}
        \int_{0}^{T} \int_{\Omega}u^{2} d x d t,  
\end{equation}
constrained to the equations:
\begin{equation}
	\begin{dcases}
		\dt {v} - \Delta v + \nabla
                p =u &\text { in }\Omega\times (0,T), \\
		-\nabla \cdot v =0 &\text { in }\Omega\times (0,T),\\
		v =0 &\text { on } \partial \Omega \times (0, T),\\
		v(0) = 0 &\text { in } \Omega\times \{0\},
	\end{dcases}
\end{equation}
where
\begin{equation}
{v}_d = \left[10(1+t)\left(1+\frac{1}{2} \cos (4 \pi t-\pi)\right ), 0\right] \in
L^2(0,T; [L^2(\Omega)]^2)
\end{equation}
is considered on the whole space-time
domain, and $\alpha = 10^{-5}$. For the space discretization we employed
$\mathbb{P}^2-\mathbb{P}^1$ polynomials for the velocity and pressure
fields both for the state and the adjoint variables, while for the
control we used $\mathbb P^2$ polynomials. Thus, the FE dimension is
$674$ for state/adjoint velocity and control variables, while it is
$337$ for the state and adjoint pressure.  
In terms of time discretization, in the interval $[0,1]$ we considered
$N_t = 50$ time instances, with $\Delta t = 0.02$. The final dimension
of the system is $134800$.  \\
In the left panel of Figure~\ref{fig:recon_pred_error_stokes} we plotted the relative
reconstruction error $E_k$ for all the variables of interest and for
every $k=1, \dots, 50$. The error at
the first time step is exactly $0$, while the spike at the last time
instant is due to the nature of DMD which divides the snapshots data
into two shifted matrices and we have the accumulation of all the
residuals. The maximum error is registered for the adjoint velocity
and it is below $1.5\%$, while on average all the variables present an
error below $0.6\%$. We remark that the adjoint pressure snapshots
have been normalized removing the mean of all the states before
applying the DMDc. Even if such preprocessing can be done only for
reconstruction, we emphasize that the adjoint pressure in this case is
not a variable of interest for prediction purposes since only the adjoint velocity affects the control in the optimality equation. Nonetheless we
present its relative reconstruction error in order to show that the partitioned approach is able to
deal also with such variable. 

In Figure~\ref{fig:recon_pred_error_stokes} (right panel) we plotted the $L^2$ mean
relative prediction error over $20$ future states for velocity,
pressure, and adjoint velocity, varying the dimension of the train
dataset. We notice an almost stationary behaviour of the mean error
for all the variables of interest, with all the values below
$0.8\%$. This suggests that the future state prediction is robust with
respect to the amount of available data.

\review{
\subsection{Speedup considerations}
The computational gain provided by the data-driven approach is
considerable for both the test problems. The one-shot approach
takes between $13$ and $18$ minutes, depending on the problem at hand,
in order to compute all the variables. The partitioned
approach exploiting DMDc, instead, requires less than $1$ second for
the construction of the reduced operators, and approximately $2$
seconds to predict all the quantities on interest for the entire time
spans considered. This results in a speedup $\sim \mathcal{O}(10^2)$.

}

\section{Conclusions and perspectives}
\label{sec:the_end}

In this work we showed the potential of data-driven methods for the
resolution of time-dependent optimal control problems with PDE
constraints.
A further improvement of the work might concern a deeper analysis of
the sensitivity of the problems with respect to the penalization
parameter $\alpha$ which drastically changes the magnitude of the
adjoint variable. 
Another possible approach could be to use the data-driven sparse
identification of nonlinear dynamics with control
method~\cite{brunton2016sparse}, instead of DMDc.
Lastly, we stress that even if the test cases are quite academic, we
believe in the potential of such an approach in more complex setting
based on data collection. Indeed, this represents a first attempt for a
deeper study of this formulation, which might be of interest in many
interdisciplinary research fields.

\section*{Acknowledgements}
This work was partially funded by European
Union Funding for Research and Innovation --- Horizon 2020 Program --- in
the framework of European Research Council Executive Agency: H2020 ERC
CoG 2015 AROMA-CFD project 681447 ``Advanced Reduced Order Methods with
Applications in Computational Fluid Dynamics'' P.I. Professor
Gianluigi Rozza.

\bibliographystyle{abbrvurl}

\begin{thebibliography}{10}

\bibitem{multiphenics}
multiphenics - easy prototyping of multiphysics problems in fenics,
  \href{https://mathlab.sissa.it/multiphenics}{https://mathlab.sissa.it/multiphenics}.

\bibitem{bai2020dynamic}
Z.~Bai, E.~Kaiser, J.~L. Proctor, J.~N. Kutz, and S.~L. Brunton.
\newblock Dynamic mode decomposition for compressive system identification.
\newblock {\em AIAA Journal}, 58(2):561--574, 2020.
\newblock \href {https://doi.org/10.2514/1.J057870}
  {\path{doi:10.2514/1.J057870}}.

\bibitem{Ballarin2017}
F.~Ballarin, E.~Faggiano, A.~Manzoni, A.~Quarteroni, G.~Rozza, S.~Ippolito,
  C.~Antona, and R.~Scrofani.
\newblock Numerical modeling of hemodynamics scenarios of patient-specific
  coronary artery bypass grafts.
\newblock {\em Biomechanics and Modeling in Mechanobiology}, 16(4):1373--1399,
  Aug 2017.
\newblock \href {https://doi.org/10.1007/s10237-017-0893-7}
  {\path{doi:10.1007/s10237-017-0893-7}}.

\bibitem{BALLARIN2022307}
F.~Ballarin, G.~Rozza, and M.~Strazzullo.
\newblock {Chapter 9 - Space-time POD-Galerkin approach for parametric flow
  control}.
\newblock In E.~Trélat and E.~Zuazua, editors, {\em Numerical Control: Part
  A}, volume~23 of {\em Handbook of Numerical Analysis}, pages 307--338.
  Elsevier, 2022.
\newblock \href {https://doi.org/10.1016/bs.hna.2021.12.009}
  {\path{doi:10.1016/bs.hna.2021.12.009}}.

\bibitem{Benzi}
M.~Benzi, G.~H. Golub, and J.~Liesen.
\newblock Numerical solution of saddle point problems.
\newblock {\em Acta Numerica}, 14:1--137, 2005.
\newblock \href {https://doi.org/10.1017/S0962492904000212}
  {\path{doi:10.1017/S0962492904000212}}.

\bibitem{brunton2019data}
S.~L. Brunton and J.~N. Kutz.
\newblock {\em {Data-Driven Science and Engineering: Machine Learning,
  Dynamical Systems, and Control}}.
\newblock Cambridge University Press, 2019.
\newblock \href {https://doi.org/10.1017/9781108380690}
  {\path{doi:10.1017/9781108380690}}.

\bibitem{brunton2016sparse}
S.~L. Brunton, J.~L. Proctor, and J.~N. Kutz.
\newblock {Sparse identification of nonlinear dynamics with control (SINDYc)}.
\newblock {\em IFAC-PapersOnLine}, 49(18):710--715, 2016.
\newblock \href {https://doi.org/10.1016/j.ifacol.2016.10.249}
  {\path{doi:10.1016/j.ifacol.2016.10.249}}.

\bibitem{dede2007optimal}
L.~Ded\`e.
\newblock Optimal flow control for {N}avier-{S}tokes equations: Drag
  minimization.
\newblock {\em International Journal for Numerical Methods in Fluids},
  55(4):347--366, 2007.
\newblock \href {https://doi.org/10.1002/fld.1464}
  {\path{doi:10.1002/fld.1464}}.

\bibitem{dede2010reduced}
L.~Ded{\`e}.
\newblock Reduced basis method and a posteriori error estimation for
  parametrized linear-quadratic optimal control problems.
\newblock {\em SIAM Journal on Scientific Computing}, 32(2):997--1019, 2010.
\newblock \href {https://doi.org/10.1137/090760453}
  {\path{doi:10.1137/090760453}}.

\bibitem{delfour2011shapes}
M.~C. Delfour and J.~Zol{\'e}sio.
\newblock {\em Shapes and geometries: metrics, analysis, differential calculus,
  and optimization}, volume~22.
\newblock SIAM, Philadelphia, 2011.
\newblock \href {https://doi.org/10.1137/1.9780898719826}
  {\path{doi:10.1137/1.9780898719826}}.

\bibitem{demo2018shape}
N.~Demo, M.~Tezzele, G.~Gustin, G.~Lavini, and G.~Rozza.
\newblock Shape optimization by means of proper orthogonal decomposition and
  dynamic mode decomposition.
\newblock In {\em Technology and Science for the Ships of the Future:
  Proceedings of NAV 2018: 19th International Conference on Ship \& Maritime
  Research}, pages 212--219. IOS Press, 2018.
\newblock \href {https://doi.org/10.3233/978-1-61499-870-9-212}
  {\path{doi:10.3233/978-1-61499-870-9-212}}.

\bibitem{demo2018pydmd}
N.~Demo, M.~Tezzele, and G.~Rozza.
\newblock {PyDMD: Python Dynamic Mode Decomposition}.
\newblock {\em The Journal of Open Source Software}, 3(22):530, 2018.
\newblock \href {https://doi.org/10.21105/joss.00530}
  {\path{doi:10.21105/joss.00530}}.

\bibitem{Fevola2021}
E.~Fevola, F.~Ballarin, L.~Jiménez-Juan, S.~Fremes, S.~Grivet-Talocia,
  G.~Rozza, and P.~Triverio.
\newblock An optimal control approach to determine resistance-type boundary
  conditions from in-vivo data for cardiovascular simulations.
\newblock {\em International Journal for Numerical Methods in Biomedical
  Engineering}, 37(10), 2021.
\newblock \href {https://doi.org/10.1002/cnm.3516}
  {\path{doi:10.1002/cnm.3516}}.

\bibitem{Carere2021261}
C.~G., S.~M., B.~F., R.~G., and S.~R.
\newblock A weighted pod-reduction approach for parametrized pde-constrained
  optimal control problems with random inputs and applications to environmental
  sciences.
\newblock {\em Computers and Mathematics with Applications}, 102:261 – 276,
  2021.
\newblock \href {https://doi.org/10.1016/j.camwa.2021.10.020}
  {\path{doi:10.1016/j.camwa.2021.10.020}}.

\bibitem{Glas2017}
S.~Glas, A.~Mayerhofer, and K.~Urban.
\newblock {\em Two Ways to Treat Time in Reduced Basis Methods}, pages 1--16.
\newblock Springer International Publishing, Cham, 2017.
\newblock \href {https://doi.org/10.1007/978-3-319-58786-8\_1}
  {\path{doi:10.1007/978-3-319-58786-8\_1}}.

\bibitem{goldschmidt2021bilinear}
A.~Goldschmidt, E.~Kaiser, J.~L. Dubois, S.~L. Brunton, and J.~N. Kutz.
\newblock Bilinear dynamic mode decomposition for quantum control.
\newblock {\em New Journal of Physics}, 23(3):033035, 2021.
\newblock \href {https://doi.org/10.1088/1367-2630/abe972}
  {\path{doi:10.1088/1367-2630/abe972}}.

\bibitem{makinen}
J.~Haslinger and R.~A.~E. M{\"a}kinen.
\newblock {\em Introduction to shape optimization: theory, approximation, and
  computation}.
\newblock Advances in Design and Control. SIAM, Philadelphia, 2003.
\newblock \href {https://doi.org/10.1137/1.9780898718690}
  {\path{doi:10.1137/1.9780898718690}}.

\bibitem{HinzeStokes}
M.~Hinze, M.~K{\"o}ster, and S.~Turek.
\newblock A hierarchical space-time solver for distributed control of the
  {S}tokes equation.
\newblock {\em Technical Report, SPP1253-16-01}, 2008.

\bibitem{HinzeNS}
M.~Hinze, M.~K{\"o}ster, and S.~Turek.
\newblock A space-time multigrid method for optimal flow control.
\newblock In {\em Constrained optimization and optimal control for partial
  differential equations}, page 147. Springer, 2012.
\newblock \href {https://doi.org/10.1007/978-3-0348-0133-1\_8}
  {\path{doi:10.1007/978-3-0348-0133-1\_8}}.

\bibitem{hinze2008optimization}
M.~Hinze, R.~Pinnau, M.~Ulbrich, and S.~Ulbrich.
\newblock {\em Optimization with {PDE} constraints}, volume~23 of {\em
  Mathematical Modelling: Theory and Applications}.
\newblock Springer Netherlands, 2009.
\newblock \href {https://doi.org/10.1007/978-1-4020-8839-1}
  {\path{doi:10.1007/978-1-4020-8839-1}}.

\bibitem{Iapichino2}
L.~Iapichino, S.~Trenz, and S.~Volkwein.
\newblock Reduced-order multiobjective optimal control of semilinear parabolic
  problems.
\newblock In B.~Karas{\"o}zen, M.~Manguo{\u{g}}lu, M.~Tezer-Sezgin,
  S.~G{\"o}ktepe, and {\"O}.~U{\u{g}}ur, editors, {\em Numerical Mathematics
  and Advanced Applications ENUMATH 2015}, pages 389--397, Cham, 2016. Springer
  International Publishing.
\newblock \href {https://doi.org/10.1007/978-3-319-39929-4\_37}
  {\path{doi:10.1007/978-3-319-39929-4\_37}}.

\bibitem{kutz2016dynamic}
J.~N. Kutz, S.~L. Brunton, B.~W. Brunton, and J.~L. Proctor.
\newblock {\em {Dynamic Mode Decomposition: Data-Driven Modeling of Complex
  Systems}}.
\newblock SIAM, 2016.
\newblock \href {https://doi.org/10.1137/1.9781611974508}
  {\path{doi:10.1137/1.9781611974508}}.

\bibitem{LassilaManzoniQuarteroniRozza2013a}
T.~Lassila, A.~Manzoni, A.~Quarteroni, and G.~Rozza.
\newblock A reduced computational and geometrical framework for inverse
  problems in hemodynamics.
\newblock {\em International Journal for Numerical Methods in Biomedical
  Engineering}, 29(7):741--776, 2013.
\newblock \href {https://doi.org/10.1002/cnm.2559}
  {\path{doi:10.1002/cnm.2559}}.

\bibitem{leugering2014trends}
G.~Leugering, P.~Benner, S.~Engell, A.~Griewank, H.~Harbrecht, M.~Hinze,
  R.~Rannacher, and S.~Ulbrich.
\newblock {\em Trends in {PDE} constrained optimization}, volume 165 of {\em
  International Series of Numerical Mathematics}.
\newblock Springer, New York, 2014.
\newblock \href {https://doi.org/10.1007/978-3-319-05083-6}
  {\path{doi:10.1007/978-3-319-05083-6}}.

\bibitem{fenics}
A.~Logg, K.~Mardal, and G.~Wells.
\newblock {\em Automated Solution of Differential Equations by the Finite
  Element Method}.
\newblock Springer-Verlag, Berlin, 2012.

\bibitem{mohammadi2010applied}
B.~Mohammadi and O.~Pironneau.
\newblock {\em Applied shape optimization for fluids}.
\newblock Oxford University Press, New York, 2010.
\newblock \href {https://doi.org/10.1093/acprof:oso/9780199546909.001.0001}
  {\path{doi:10.1093/acprof:oso/9780199546909.001.0001}}.

\bibitem{narasingam2017development}
A.~Narasingam and J.~S.-I. Kwon.
\newblock Development of local dynamic mode decomposition with control:
  Application to model predictive control of hydraulic fracturing.
\newblock {\em Computers \& Chemical Engineering}, 106:501--511, 2017.
\newblock \href {https://doi.org/10.1016/j.compchemeng.2017.07.002}
  {\path{doi:10.1016/j.compchemeng.2017.07.002}}.

\bibitem{negri2015reduced}
F.~Negri, A.~Manzoni, and G.~Rozza.
\newblock Reduced basis approximation of parametrized optimal flow control
  problems for the {S}tokes equations.
\newblock {\em Computers \& Mathematics with Applications}, 69(4):319--336,
  2015.
\newblock \href {https://doi.org/10.1016/j.camwa.2014.12.010}
  {\path{doi:10.1016/j.camwa.2014.12.010}}.

\bibitem{proctor2016dynamic}
J.~L. Proctor, S.~L. Brunton, and J.~N. Kutz.
\newblock {Dynamic Mode Decomposition with Control}.
\newblock {\em SIAM Journal on Applied Dynamical Systems}, 15(1):142--161,
  2016.
\newblock \href {https://doi.org/10.1137/15M1013857}
  {\path{doi:10.1137/15M1013857}}.

\bibitem{quarteroni2005numerical}
A.~Quarteroni, G.~Rozza, L.~Ded{\`e}, and A.~Quaini.
\newblock Numerical approximation of a control problem for advection-diffusion
  processes.
\newblock In {\em Ceragioli F., Dontchev A., Futura H., Marti K., Pandolfi L.
  (eds) System Modeling and Optimization. International Federation for
  Information Processing, CSMO Conference on System Modeling and Optimization},
  pages vol 199, 261--273. Springer, Boston, 2005.

\bibitem{quarteroni2007reduced}
A.~Quarteroni, G.~Rozza, and A.~Quaini.
\newblock Reduced basis methods for optimal control of advection-diffusion
  problems.
\newblock In {\em Advances in Numerical Mathematics}, number 2006-003 in
  CMCS-CONF, pages 193--216. RAS and University of Houston, 2007.

\bibitem{morhandbook2020}
G.~Rozza, M.~Hess, G.~Stabile, M.~Tezzele, and F.~Ballarin.
\newblock {Basic Ideas and Tools for Projection-Based Model Reduction of
  Parametric Partial Differential Equations}.
\newblock In P.~Benner, S.~Grivet-Talocia, A.~Quarteroni, G.~Rozza, W.~H.~A.
  Schilders, and L.~M. Silveira, editors, {\em Model Order Reduction},
  volume~2, chapter~1, pages 1--47. De Gruyter, Berlin, Boston, 2020.
\newblock \href {https://doi.org/10.1515/9783110671490-001}
  {\path{doi:10.1515/9783110671490-001}}.

\bibitem{schoberl2007symmetric}
J.~Sch{\"o}berl and W.~Zulehner.
\newblock Symmetric indefinite preconditioners for saddle point problems with
  applications to {PDE}-constrained optimisation problems.
\newblock {\em SIAM Journal on Matrix Analysis and Applications},
  29(3):752--773, 2007.
\newblock \href {https://doi.org/10.1137/060660977}
  {\path{doi:10.1137/060660977}}.

\bibitem{seymen2014distributed}
Z.~K. Seymen, H.~Y{\"u}cel, and B.~Karas{\"o}zen.
\newblock Distributed optimal control of time-dependent
  diffusion--convection--reaction equations using space--time discretization.
\newblock {\em Journal of Computational and Applied Mathematics}, 261:146--157,
  2014.
\newblock \href {https://doi.org/10.1016/j.cam.2013.11.006}
  {\path{doi:10.1016/j.cam.2013.11.006}}.

\bibitem{Stoll}
M.~Stoll and A.~Wathen.
\newblock All-at-once solution of time-dependent {S}tokes control.
\newblock {\em J. Comput. Phys.}, 232(1):498--515, Jan. 2013.
\newblock \href {https://doi.org/10.1016/j.jcp.2012.08.039}
  {\path{doi:10.1016/j.jcp.2012.08.039}}.

\bibitem{Strazzullo1}
M.~Strazzullo, F.~Ballarin, R.~Mosetti, and G.~Rozza.
\newblock Model reduction for parametrized optimal control problems in
  environmental marine sciences and engineering.
\newblock {\em SIAM Journal on Scientific Computing}, 40(4):B1055--B1079, 2018.
\newblock \href {https://doi.org/10.1137/17M1150591}
  {\path{doi:10.1137/17M1150591}}.

\bibitem{Strazzullo2}
M.~Strazzullo, F.~Ballarin, and G.~Rozza.
\newblock {{POD}-{G}alerkin Model Order Reduction for Parametrized Time
  Dependent Linear Quadratic Optimal Control Problems in Saddle Point
  Formulation}.
\newblock {\em Journal of Scientific Computing}, 83(55), 2020.
\newblock \href {https://doi.org/10.1007/s10915-020-01232-x}
  {\path{doi:10.1007/s10915-020-01232-x}}.

\bibitem{StrazzulloRB}
M.~Strazzullo, F.~Ballarin, and G.~Rozza.
\newblock A {C}ertified {R}educed {B}asis method for linear parametrized
  parabolic optimal control problems in space-time formulation.
\newblock Submitted, 2021,
  \href{https://arxiv.org/abs/2103.00460}{https://arxiv.org/abs/2103.00460}.

\bibitem{StrazzulloSWE}
M.~Strazzullo, F.~Ballarin, and G.~Rozza.
\newblock {POD-Galerkin model order reduction for parametrized nonlinear
  time-dependent optimal flow control: an application to shallow water
  equations}.
\newblock {\em Journal of Numerical Mathematics}, 30(1):63--84, 2022.
\newblock \href {https://doi.org/10.1515/jnma-2020-0098}
  {\path{doi:10.1515/jnma-2020-0098}}.

\bibitem{ZakiaMaria}
M.~Strazzullo, Z.~Zainib, F.~Ballarin, and G.~Rozza.
\newblock Reduced order methods for parametrized nonlinear and time dependent
  optimal flow control problems: towards applications in biomedical and
  environmental sciences.
\newblock {\em Numerical Mathematics and Advanced Applications ENUMATH 2019},
  2021.
\newblock \href {https://doi.org/10.1007/978-3-030-55874-1_83}
  {\path{doi:10.1007/978-3-030-55874-1_83}}.

\bibitem{tezzele2018ecmi}
M.~Tezzele, N.~Demo, A.~Mola, and G.~Rozza.
\newblock {An integrated data-driven computational pipeline with model order
  reduction for industrial and applied mathematics}.
\newblock In M.~G{\"u}nther and W.~Schilders, editors, {\em {Novel Mathematics
  Inspired by Industrial Challenges}}, number~X in Mathematics in Industry.
  Springer International Publishing, 2022.

\bibitem{tezzele2020enhancing}
M.~Tezzele, N.~Demo, G.~Stabile, A.~Mola, and G.~Rozza.
\newblock {Enhancing CFD predictions in shape design problems by model and
  parameter space reduction}.
\newblock {\em Advanced Modeling and Simulation in Engineering Sciences},
  7(40), 2020.
\newblock \href {https://doi.org/10.1186/s40323-020-00177-y}
  {\path{doi:10.1186/s40323-020-00177-y}}.

\bibitem{troltzsch2010optimal}
F.~Tr{\"o}ltzsch.
\newblock {Optimal Control of Partial Differential Equations: Theory, Methods
  and Applications}.
\newblock {\em Graduate Studies in Mathematics}, 112, 2010.
\newblock \href {https://doi.org/10.1090/gsm/112} {\path{doi:10.1090/gsm/112}}.

\bibitem{urban2012new}
K.~Urban and A.~T. Patera.
\newblock A new error bound for reduced basis approximation of parabolic
  partial differential equations.
\newblock {\em Comptes Rendus Mathematique}, 350(3-4):203--207, 2012.
\newblock \href {https://doi.org/10.1016/j.crma.2012.01.026}
  {\path{doi:10.1016/j.crma.2012.01.026}}.

\bibitem{yano2014space}
M.~Yano.
\newblock A space-time {P}etrov--{G}alerkin certified reduced basis method:
  Application to the {B}oussinesq equations.
\newblock {\em SIAM Journal on Scientific Computing}, 36(1):A232--A266, 2014.
\newblock \href {https://doi.org/10.1137/120903300}
  {\path{doi:10.1137/120903300}}.

\bibitem{yano2014space1}
M.~Yano, A.~T. Patera, and K.~Urban.
\newblock A space-time hp-interpolation-based certified reduced basis method
  for {B}urgers' equation.
\newblock {\em Mathematical Models and Methods in Applied Sciences},
  24(09):1903--1935, 2014.
\newblock \href {https://doi.org/10.1142/S0218202514500110}
  {\path{doi:10.1142/S0218202514500110}}.

\bibitem{Zakia}
Z.~Zainib, F.~Ballarin, S.~Fremes, P.~Triverio, L.~Jim\'{e}nez-Juan, and
  G.~Rozza.
\newblock Reduced order methods for parametric optimal flow control in coronary
  bypass grafts, toward patient-specific data assimilation.
\newblock {\em International Journal for Numerical Methods in Biomedical
  Engineering}, page e3367, 2020.
\newblock \href {https://doi.org/10.1002/cnm.3367}
  {\path{doi:10.1002/cnm.3367}}.

\end{thebibliography}

\end{document}